\documentclass[a4paper, 12pt]{article}
\usepackage{amsmath, amsthm, mathrsfs, amssymb, verbatim, bbm, url}
  \usepackage{amscd}
\usepackage{latexsym,amssymb}
\usepackage{graphicx}
\usepackage{color}
\usepackage{enumerate}

    \newcommand{\ix}{\mathcal{X}}
      \newcommand{\iy}{\mathcal{Y}}

\newcommand{\oxvir}[0]{\mathcal{O}_{\mathcal{X}}^{vir}}
\newcommand{\oxmvir}[0]{\mathcal{O}_{\mathcal{X}_\mu}^{vir}}
\newtheorem{virkaw}{Lemma}[section]
\newtheorem{virkaw22}[virkaw]{Theorem}
\newtheorem{arem0}[virkaw]{Remark}
\newtheorem{cones}{Lemma}[section]
\newtheorem{zerolema}[cones]{Lemma}
\newtheorem{cones11}[cones]{Remark}
\newtheorem{arem1}[virkaw]{Remark}
\newtheorem{arem2}[virkaw]{Remark}

\begin{document}

\begin{center}
\Large{ A virtual Kawasaki formula} \vspace{1 cm}

\large{Valentin Tonita} 
\end{center}
\begin{abstract}
 Kawasaki's formula is a tool to compute holomorphic Euler characteristics of vector bundles on a compact orbifold $\ix$. Let $\ix$ be an orbispace with perfect obstruction theory which admits an embedding in a smooth orbifold. One can then construct the virtual structure sheaf and the virtual fundamental class of $\ix$. In this paper we prove that Kawasaki's formula ``behaves well'' with working  ``virtually'' on $\ix$ in the following sense: if we replace the structure sheaves, tangent and normal bundles in the formula by their virtual counterparts then Kawasaki's formula stays true. Our motivation comes from studying the quantum K-theory of a complex manifold $X$ (see \cite{gito}), with the formula applied to Kontsevich' moduli spaces of genus $0$ stable maps to $X$.
 
 \end{abstract}
 
 \section{Introduction}
 Given a manifold $\ix$ and a vector bundle $V$ on $\ix$ then Hirzebruch-Riemann-Roch formula states that:
  \begin{align*}
 \chi(\ix, V) =\int_\ix ch(V)Td(T_\ix ) . 
  \end{align*}
  
   In \cite{kaw} Kawasaki generalized this formula to the case when $\ix$ is an orbifold. He reduces the computation of Euler characteristics on $\ix$ to computation of certain cohomological integrals on \emph{the inertia orbifold} $I\ix$: 
    \begin{align}
    \chi(\ix, V) = \sum_\mu \frac{1}{m_\mu}\int_{\ix_\mu} Td(T_{\ix_\mu})ch \left(\frac{Tr(V)}{Tr(\Lambda^\bullet N^*_\mu)}\right).\label{001} 
    \end{align}
 We explain below the ingredients in the formula:  \\
 $I\ix$ is defined as follows: around any point $p\in \ix$ there is a local chart $(\widetilde{U}_p, G_p)$ such that locally $\ix$ is represented as the quotient of $\widetilde{U}_p$ by $G_p$. Consider the set of conjugacy classes $(1)=(h_p^1)$, $(h_p^2)$, $\ldots$, $(h_p^{n_p})$ in $G_p$. Define:
      \begin{align*}
      I\ix:=\{\left(p,(h_p^i)\right) \quad \vert \quad i=1,2,\ldots, n_p \}.
       \end{align*}
   Pick an element $h_p^i$ in each conjugacy class. Then a local chart on $I\ix$ is given by:
    \begin{align*}
       \coprod_{i=1}^{n_p} \widetilde{U}_p^{(h_p^i)}/ Z_{G_p}(h_p^i),
     \end{align*}
     where $Z_{G_p}(h_p^i)$ is the centralizer of $h_p^i$ in $G_p$.
   Denote by $\ix_\mu$ the connected components of the inertia orbifold (we'll often refer to them as Kawasaki strata). The multiplicity $m_\mu$ associated to each $\ix_\mu$ is given by:
      \begin{align*}
     m_\mu:= \left\vert ker\left(Z_{G_p}(g)\rightarrow Aut(\widetilde{U}_p^g)\right) \right\vert . 
     \end{align*} 
   
   For a vector bundle $V$ we will denote by $V^*$ the dual bundle to $V$. The restriction of $V$ to $\ix_\mu$ decomposes in characters of the $g$ action. Let $E_r^{(l)}$ be the subbundle of the restriction of $E$ to $\ix_\mu$ on which $g$ acts with eigenvalue $e^{\frac{2\pi i l}{r}}$. Then the  trace $Tr(V)$ is defined to be the orbibundle whose fiber over the point $(p, (g))$ of $\ix_\mu$ is :
    \begin{align*}
 Tr(V):= \sum_{l} e^{\frac{2\pi i l}{r}} E^{(l)}_r . 
    \end{align*} 
 Finally, $\Lambda^\bullet N^*_\mu$ is the K-theoretic Euler class of the normal bundle $N_\mu$ of $\ix_\mu$ in $\ix$. $Tr(\Lambda^\bullet N^*_\mu)$ is invertible because the symmetry $g$ acts with eigenvalues different from $1$ on the normal bundle to the fixed point locus.     
  We call the terms corresponding to the identity component in the formula \emph{fake Euler characteristics}:
     \begin{align*}
  \chi^f(\ix,V)= \int_\ix ch(V) Td(T_\ix) .      
     \end{align*}
   In the case where $\ix$ is a global quotient formula $(\ref{001})$ is the Lefschetz fixed point formula. 
  
  Now let $\ix$ be a compact, complex orbispace (Deligne-Mumford stack) with a perfect obstruction theory $E^{-1}\to E^{0}$. This gives rise to the intrinsic normal cone, which is embedded in $E_1$ - the dual bundle to $E^{-1}$ (see \cite{liti}, also \cite{behfan}). The virtual structure sheaf $\oxvir$ was defined in \cite{ypl2}  as the K-theoretic pull-back by the zero section of the structure sheaf of this cone. Let $I\ix =\coprod_\mu \ix_\mu$ be the inertia orbifold of $\ix$. We denote by $i_\mu$ the inclusion of a stratum $\ix_\mu$ in $\ix$. For a bundle $V$ on $\ix$ we write $i_\mu^* V =V_\mu^f\oplus V_\mu^m$ for its decomposition as the direct sum  of the fixed part and the moving part under the action of the symmetry associated to $\ix_\mu$. To avoid ugly notation we will often simply write $V^m, V^f$. The virtual normal bundle to $\ix_\mu$ in $\ix$ is defined as $[E_{0}^m]-[E_{1}^m]$. We will in addition assume that $\ix$ admits an embedding $j$ in a smooth compact orbifold $\iy$. This is always true for the moduli spaces of genus $0$ stable maps $X_{0,n,d}$ because an embedding $X\hookrightarrow \mathbb{P}^N$ induces an embedding $X_{0,n,d}\hookrightarrow (\mathbb{P}^N)_{0,n,d}$.

    \begin{virkaw22} \label{mmmm}
  {\em Denote by $N^{vir}_\mu$ the virtual normal bundle of $\ix_\mu$ in $\ix$. Then } 
  \begin{align}
  \chi\left(\ix, j^* (V)\otimes\oxvir\right) = \sum_\mu \frac{1}{m_\mu} \chi^f\left(\ix_\mu, \frac{Tr(V_\mu\otimes\oxmvir)}{Tr\left(\Lambda^\bullet (N^{vir}_\mu)^*\right)} \right).
  \end{align}
   \end{virkaw22}
  \begin{arem0}
 {\em A perfect obstruction theory $E^{-1}\to E^0$ on $\ix$ induces canonically a perfect obstruction theory on $\ix_\mu$ by taking the fixed part of the complex $E_\mu^{-1,f}\to E_\mu^{0,f}$. The proof is the same as that of Proposition $1$ in \cite{grpa}. This is then used to define the sheaf $\oxmvir$.} 
  \end{arem0} 
  
 \begin{arem2}
 {\em It is proved in \cite{fago} that if $\ix$ is a scheme, the  Grothendieck-Riemann-Roch theorem is compatible with  virtual fundamental classes and virtual fundamental sheaves i.e.: }
 \begin{align*}
 \chi^f(\ix, V\otimes \oxvir) = \int_{[\ix]} ch(V\otimes \oxvir)\cdot Td(T^{vir}) 
  \end{align*}
 {\em where $[\ix]$ is the virtual fundamental class of $\ix$ and $T^{vir}$ is its virtual tangent bundle. Their arguments carry over to the case when $\ix$ is a stack.} 
 \end{arem2}
  \begin{arem1}
 {\em The bundles $V$ to which we apply  Theorem $\ref{mmmm}$ in \cite{gito} are (sums and products of) cotangent line bundles $L_i$ and evaluation classes $ev_i^*(a_i)$. They are pull-backs of the corresponding bundles on $(\mathbb{P}^N)_{0,n,d}$.}
  \end{arem1}

$\mathbf{Acknowledgements}$. I would like to thank Alexander Givental for suggesting the problem and for useful discussions. Thanks are also due to Yuan-Pin Lee who patiently answered my questions on the material in his work \cite{ypl2} and to Hsian-Hua Tseng who read a preliminary draft of the paper.  
 \section{Proof of Theorem $\ref{mmmm}$} 
 
  
  Before proving Theorem $\ref{mmmm}$ we recall a couple of background facts and lemmata on K-theory which we will use.
  
 Let $K_0(X)$ be the Grothendieck group of coherent sheaves on $X$. Given a map $f:X\to Y$, the K-theoretic pullback $f^*(\mathcal{F}):K_0(Y)\to K_0(X)$ is defined as the alternating sum of derived functors $Tor_{\mathcal{O}_Y}^i(\mathcal{F},\mathcal{O}_X)$, provided that the sum is finite. This is always true for instance if $f$ is flat or if it is a regular embedding. 
 
For any fiber square:
     \begin{displaymath}
     \begin{CD}
     V' @>>>V\\
      @VVV               @VVV                          \\
 B'     @> i >> B
             \end{CD}
 \end{displaymath}
 with $i$ a regular embedding one can define K-theoretic refined Gysin homomorphisms $i^! : K_0(V)\to K_0(V')$ (see \cite{ypl2}). One way to define the map $i^!$ is the following: the class $i_*(\mathcal{O}_{B'})\in K^0(B)$ has a finite resolution of vector bundles, which is exact off $B'$. We pull it back to $V$ and then cap (i.e. tensor product) with classes in $K_0(V)$, to get a class on $K_0(V)$ with homology supported on $V'$, which we can regard as an element of $K_0(V')$, because there is a canonical isomorphism between complexes on $V$ with homology supported on $V'$  and $K_0(V')$. 

 In the following two lemmata $X,Y,Y'$ are assumed DM stacks. We will use the following result:
\begin{cones} \label{le0001}

  Consider the diagram:
  
     \begin{displaymath}
     \begin{CD}
     \iota^*C_{X/Y} @>>>  C_{X/Y} \\
      @VVV                @VVV      \\
     X' @> \iota >>X\\
      @VVV               @V j VV                          \\
  Y'     @> i >> Y
             \end{CD}
 \end{displaymath}
  with $i$ a regular embedding and $j$ an embedding, $C_{X/Y}$ is the normal cone of $X$ in $Y$ and both squares are fiber diagrams. Then:
    \begin{align}  
    i^{!} [\mathcal{O}_{C_{X/Y}}] = [\mathcal{O}_{C_{X'/Y'}}] \in K_0( \iota^*C_{X/Y}) \label{con1}.
    \end{align} 
   \end{cones} 

 This is stated and proved in \cite{ypl2} (Lemma 2). The proof is based on a more general statement (Lemma 1 of \cite{ypl2}), which has been worked  out in \cite{krs} on the level of Chow rings. Since K-theoretic statements are stronger, we give below the key-ingredient which allows one to carry over Kresch's proof to K-theory:
 
  \begin{zerolema} \label{le0000}
  Let $f: X \to Y$ be a closed embedding let $g:Y\to \mathbb{P}^1 $ be a surjection such that $g\circ f$ is flat. Denote by $X_0$ and $Y_0$ the fibers over $0$ of $g\circ f$ and $g$ respectively. Moreover assume that the restriction of $f$ to $X\setminus X_0$ is an isomorphism. Then if $i$ is the inclusion of $\{0\}$ in $\mathbb{P}^1$, $i^!(\mathcal{O}_{Y}) =\mathcal{O}_{X_0}\in K_0(Y_0).$
   \end{zerolema}
   
  Proof: the skyscraper sheaves at all points of $\mathbb{P}^1$ represent the same element in $K_0(\mathbb{P}^1)$ , hence if we pull-back a resolution of any point $P\in \mathbb{P}^1$ by $g$  we get the same elements of $K_0(Y)$. On the other hand since $f$ is an isomorphism above $\mathbb{P}^1\setminus \{0\}$, pulling-back by $g$ of the structure sheaf of a point $P\neq 0$ is the same as pulling back by $g\circ f$ followed by $f_*$. By what we said above we can replace $P$ with $0$.  Now from the flatness of $g\circ f$ above $0$ the pull-back of the structure sheaf of $0$ by $g\circ f$ is the structure sheaf of the fiber $X_0$. The result then follows from  the definition of $i^!$.  
    
  \begin{cones11}
 {\em  Lemma $\ref{le0000}$ allows one to show  Lemma $\ref{le0001}$: intermediately one shows, following \cite{krs},(notation is as in Lemma $\ref{le0001}$) that $[\mathcal{O}_{C_1}]=[\mathcal{O}_{C_2}]$ in $K_0(C_{X'}Y\times_Y C_{X}Y )$ where $C_1:= C_{i^* C_{X}Y}(C_{X}Y)$ and $C_2:= C_{j^* C_{Y'}Y}(C_{Y'}Y)$.}
  
  \end{cones11}



We now go on to prove Theorem $\ref{mmmm}$. We have:
    \begin{align*}
   \chi(\ix, j^* V\otimes\oxvir) = \chi(\iy, V\otimes j_*\oxvir ).  
    \end{align*}
 Kawasaki's formula applied to the sheaf $V\otimes j_*\oxvir$ on $\iy$ gives:
  \begin{align}
 \chi(\iy, V\otimes j_*\oxvir ) = \sum_\mu \frac{1}{m_\mu}\chi^f\left( \iy_\mu, \frac{Tr(V_\mu \otimes i_\mu^* j_*\oxvir)}{Tr(\Lambda^\bullet N^*_\mu )} \right). \label{rel11}
 \end{align}
 From the fiber diagram:
   \begin{displaymath}
     \begin{CD}
     \ix_\mu @> i'_\mu >> \ix\\
      @V j' VV               @V j VV                          \\
 \iy_\mu     @> i_\mu >> \iy
             \end{CD}
 \end{displaymath}
     and Theorem $\mathbf{6.2}$ in  \cite{fulton} (where this is proved for Chow rings) we have $i_\mu^*j_* \oxvir = j'_* i_\mu^! \oxvir$. Plugging this in $(\ref{rel11})$ gives:
   \begin{align}
 \chi^f\left( \iy_\mu, \frac{Tr\left(V_\mu\otimes i_\mu^*j_*\oxvir\right)}{Tr(\Lambda^\bullet N^*_\mu )} \right) = \chi^f\left( \iy_\mu, \frac{Tr\left(V_\mu\otimes   j'_*i_\mu^!\oxvir\right)}{Tr(\Lambda^\bullet N^*_\mu)} \right). \label{a12}
   \end{align} 
    Let $G_\mu$ be the cyclic group generated by one element of the conjugacy class associated to $\ix_\mu$. Then we will show that: 
   \begin{align}
   Tr\left( \frac{i_\mu^!\oxvir}{\Lambda^\bullet (N^*_\mu)}\right) = Tr\left( \frac{\oxmvir }{\Lambda^\bullet (N_\mu^{vir})^*}\right) \label{a11}
      \end{align} 
    in the $G_\mu$-equivariant K-ring of $\ix_\mu$. This is essentially the computation of Section 3 in \cite{grpa} carried out in $\mathbb{C}^*$-equivariant K-theory. Relation $(\ref{a11})$ then follows by embedding the group $G_\mu$ in the torus and specializing the value of the variable $t$ in the ground ring of $\mathbb{C}^*$-equivariant K-theory to a $\vert G_\mu\vert$-root of unity.
    
   If we define a cone $D:= C_{\ix/\iy}\times_\ix E_0$, then this is a $T\iy$ cone (see \cite{behfan}). The virtual normal cone $D^{vir}$ is defined as $D/T\iy$ and $\oxvir$ is the pull-back by the zero section of the structure sheaf of $D^{vir}$. Alternatively there is a fiber diagram:
     
    \begin{displaymath}
     \begin{CD}
     T\iy @>>> D\\
      @VVV               @VVV                          \\
 \ix     @> 0_{E_1} >> E_1
             \end{CD}
 \end{displaymath}
 whre the bottom map is the zero section of $E_1$. Then one can define $\oxvir$ as $0^*_{T\iy}0_{E_1}^![\mathcal{O}_{D}]$.  We'll prove formula $(\ref{a11})$ following closely the calculation in \cite{grpa}. First by definition of $\oxvir$ and by commutativity of Gysin maps we have :
  \begin{align}
  i_\mu^! \oxvir  = i_\mu^! 0^*_{T\iy}0_{E_1}^![\mathcal{O}_{D}]
                 = 0^*_{T\iy}0_{E_1}^! i^!_\mu [\mathcal{O}_{D}].
                   \label{4a1}
  \end{align}
   We pull-back relation  $(\ref{con1})$ to $(i'_\mu)^*D = (i'_\mu)^* (C_{\ix/\iy}\times E_0)$ to get:
   \begin{align}
  i^!_\mu [\mathcal{O}_{D}] = [\mathcal{O}_{D_\mu}\times (E_0^m)^*].\label{444a1}
   \end{align}
   In the equality above we have used the fact that $D_\mu = C_{\ix_\mu/\iy_\mu}\times E_0^f$ and we identified the sheaf of sections of the bundle $E_0^m$ with the dual bundle $(E_0^m)^*$.
   Plugging $(\ref{444a1})$ in $(\ref{4a1})$ we get:
   \begin{align}
  i_\mu^! \oxvir = 0^*_{T\iy}0_{E_1}^! [\mathcal{O}_{D_\mu}\times (E_0^m)^*].\label{5555a1}
   \end{align}
   Notice that the action of $T\iy_\mu$ leaves $D_\mu\times (E_0^m)^* $ invariant (it acts trivially on $(E_0^m)^*$). Now we can write $0^*_{T\iy}=0_{T\iy_\mu^f}^*\times 0^*_{T\iy_\mu^m}$ and since $D_\mu^{vir} = D_\mu/T\iy_\mu$ we rewrite $(\ref{5555a1})$ as:
    \begin{align}
  i_\mu^! \oxvir = 0^*_{T\iy_\mu^m}0^!_{E_1}[\mathcal{O}_{D^{vir}_\mu}\times (E_0^m)^*]. \label{555a1}
    \end{align}
    
   The proof of Lemma 1 in \cite{grpa} works in our set-up as well: it uses excess intersection formula which holds in K-theory. It shows that the following relation holds in the $\mathbb{C}^*$-equivariant K-ring of $\ix_\mu$:
    \begin{align}
0^*_{T\iy_\mu^m}0^!_{E_1}[\mathcal{O}_{D^{vir}_\mu}\times (E_0^m)^*] = 0_{E_0^m}^*\left( 0_{E_1}^![\mathcal{O}_{D_\mu^{vir}}\times (E_0^m)^*]\right)\cdot \frac{\Lambda^\bullet (T\iy^m)^*}{\Lambda^\bullet(E_0^m)^*}. \label{5a1}  
    \end{align}  
    The class $0_{E_1}^![\mathcal{O}_{D_\mu^{vir}}\times E_0^m]$ lives in the $\mathbb{C}^*$-equivariant K-ring of $E^m_{0}$. The class doesn't depend on the bundle map 
    $E^m_0\to E_1^m$ so we can assume this map to be $0$. Then by excess intersection formula and the definition of $\oxmvir$ we get :
     \begin{align}
   0_{E_0^m}^*\left( 0_{E_1}^![\mathcal{O}_{D_\mu^{vir}}\times (E_0^m)^*]\right) =\oxmvir \cdot \Lambda^\bullet(E_1^m)^*.  \label{6a1}   
     \end{align}
  Formula $(\ref{6a1})$ holds because $D_\mu^{vir}\times (E_0^m)\subset E_1^f \times E_0^m$ and $0_{E_1}^!$ acts as $0^!_{E_1^f}\times 0^!_{E_1^m}$  on factors.  $0^!_{E_1^f}[\mathcal{O}_{D_\mu^{vir}}] =\oxmvir$ by definition of $\oxmvir$. By excess intersection formula applied to the fiber square:
  \begin{displaymath}
     \begin{CD}
     E_0^m  @>>>  E_0^m\\
      @V \pi VV               @VVV                          \\
 \ix_\mu     @> 0_{E^m_1} >> E^m_1
             \end{CD}
 \end{displaymath}
 we have $0_{E_0^m}^*0^!_{E_1^m}[(E_0^m)^*]= 0_{E_0^m}^*\pi^* \Lambda^\bullet(E_1^m)^* = \Lambda^\bullet(E_1^m)^*$.
     Plugging formula $(\ref{6a1})$ in $(\ref{5a1})$ (note that $N_\mu = T\iy^m_\mu$ and $N_\mu^{vir} = [E^m_{0}]-[E_{1}^m]$) and taking traces proves $(\ref{a11})$. We now  plug $(\ref{a11})$  in $(\ref{a12})$ and then pull-back to $\ix_\mu$ to get:
    \begin{align}
  \chi^f\left( \iy_\mu, \frac{Tr (V_\mu\otimes j_*i_\mu^*\oxvir)}{Tr(\Lambda^\bullet N^*_\mu)} \right) & = \chi^f \left(\iy_\mu, Tr( V_\mu)\otimes j'_*\frac{Tr(\oxmvir) }{Tr(\Lambda^\bullet (N_\mu^{vir})^*)}\right) = \nonumber \\
    & = \chi^f\left(\ix_\mu, \frac{Tr(V_\mu\otimes \oxmvir)}{Tr(\Lambda^\bullet (N^{vir}_\mu)^*)}\right) . 
    \end{align}
    
   This concludes the proof of the proposition.

\end{document}